\documentclass[11pt]{article}

\usepackage{amsmath}
\usepackage{amsthm}
\usepackage{amsfonts}
\usepackage[dvips]{graphics}
\usepackage{epsfig}
\usepackage{color}

\newtheorem{thm}{Theorem}

\newtheorem{cor}[thm]{Corollary}

\newtheorem{ex}{Application}[section]
\theoremstyle{definition} 
\title{Maximum Run Length \\ In a Toroidal Grid Graph}
\author{Margaret I. Doig \\ University of Notre Dame \\
mdoig@nd.edu\footnote{Please send all paper correspondence to:
Margaret Doig, 137 Lewis Hall, Notre Dame, IN 46556-5616}}

\begin{document}

\maketitle

\begin{abstract}

A \emph{grid graph} is a Cartesian product $G=\gamma_1 \times
\gamma_2 \times \cdots \times \gamma_k$ where the $\gamma_i$ are
cycles or paths.  The \emph{run length} of a Hamiltonian cycle in
a grid graph is defined to be the maximum number $r$ such that any
$r$ consecutive edges include no more than one edge in any
dimension. In this paper, we present several methods for producing
cycles of high run length from cycles of lower run length; sample
applications include showing that the maximum run length of $G$ is
(i) less than $k$ if each $\gamma_i$ is directed; (ii) at least
$\lfloor k/3 \rfloor +1$ if each $| \gamma_i |$ has the form
$p^{r_i}$ for some fixed prime $p$; and (iii) at least $\lfloor
k/2 \rfloor +1$ if every $| \gamma_i |$ is of the form
$2^{q_i}p^{r_i}$ for fixed prime $p$.
\end{abstract}

\section{Introduction}

We investigate the run length of a cycle or path $H$ which occurs
in a grid graph $G$.  A \emph{grid graph} is the Cartesian product
$G=\gamma_1 \times \gamma_2 \times \cdots \times \gamma_k$ where
the $\gamma_i$ are either cycles or paths.  The grid graph
$\gamma_1 \times \gamma_2 \times \cdots \times \gamma_k$ is said
to be \emph{k-dimensional}, and we sometimes refer to $\gamma_i$
as a \emph{dimension} of $G$.  The \emph{size} of the dimension is
the number of vertices in $\gamma_i$.  A \emph{torus} or
\emph{toroidal grid graph} is a grid graph where each $\gamma_i$
is a cycle. A grid graph is \emph{directed} if each $\gamma_i$ is
directed and \emph{undirected} if no $\gamma_i$ is directed.

The \emph{run length} of a cycle or path $H$ in a grid graph,
denoted $rl(H)$, is the maximum number $r$ such that any sequence
of $r$ consecutive edges in $H$ contains no more than one edge in
any one dimension. The \emph{maximum run length} of a grid graph
or $mrl(G)$ is the maximum $rl(H)$ for all Hamiltonian cycles $H$
in $G$, that is, for all cycles which pass through each vertex
exactly once.

Run length first arose in computing applications. Since every
Hamiltonian path in a torus corresponds to a listing of the
coordinates of consecutive vertices, a path on a $k$-cube
corresponds to a Gray code, i.e., a listing of all $k$-bit binary
words in such a way that any two consecutive words differ in only
one bit position. The \emph{gap} of a Gray code is the run length
of the corresponding Hamiltonian path.  Electronic
position-to-digital converters such as photon detectors use a
Hamiltonian path to match each received photon to a binary code
word representing position, and machine error decreases as the
code's gap increases. (See Goddyn et al. \cite{optrun} for a more
detailed explanation.)

The run length of Hamiltonian cycles on $k$-cubes was addressed in
1988 by Goddyn, Lawrence, and Nemeth \cite{optrun}, who
investigated the gap of $k$-bit Gray codes and gave methods for
producing them which achieve a gap that approaches $2k/3$ for
large $k$.  In 2001, Goddyn and Gvozdjak \cite{bitrun} placed a
lower bound on the gap of a Gray code which is also a lower bound
on the maximum run length of a $k$-cube $Q_k$. They showed that
$mrl(Q_k)/k \rightarrow 1$ as $k \rightarrow \infty$. In addition,
Ruskey and Sawada \cite{bent} thoroughly investigated the case of
\emph{bent} Hamiltonian paths and cycles, those whose run length
is at least 2. If $G = \gamma_1 \times \cdots \times \gamma_k$ is
a $k$-dimensional torus, then $mrl(G) \geq 2$ either when $k \geq
3$ or when $k=2$ and the sizes of the dimensions are even.  If $G$
is a $k$-dimensional non-toroidal grid graph and $k \geq 3$, then
$mrl(G) \geq 2$ if and only if the size of some dimension is even.
Ruskey and Sawada also give conditions for the existence of bent
Hamiltonian paths in graphs without bent cycles.

In this paper, we give lower bounds for $rl(H)$ and $mrl(G)$,
where $H$ is a Hamiltonian cycle in a torus $G$, by presenting
several methods for producing cycles of high run length. As
explained in Section~\ref{secdecomp}, one method considers the
Cartesian product $G$ of tori $G_1$ and $G_2$. A Hamiltonian cycle
$H$ can be constructed in $G$ from Hamiltonian cycles $H_1$ in
$G_1$ and $H_2$ in $G_2$ where $rl(H)$ is greater than either
$rl(H_1)$ or $rl(H_2)$, and we give conditions for producing such
an $H$ with run length greater than 2. These methods may be easily
applied inductively to the cases where $G= \gamma_1 \times
\gamma_2 \times \cdots \times \gamma_k$ and
$|\gamma_i|=2^{q_i}p^{r_i}$ for a fixed odd prime $p$ and for
nonnegative integers $r_i$ and $q_i$.  For situations where this
method is not applicable, Section~\ref{secdim} provides a simpler
method of combining two grid graphs $H_1$ and $H_2$ where $rl(H)$
is merely equal to $rl(H_1)$ or $rl(H_2)$. Throughout, except for
a few cases we have noted in the text, every result for finding
Hamiltonian cycles can be easily adapted to find Hamiltonian paths
as well.

\section{Composition of cycles} \label{secdecomp}

\begin{figure}
\begin{center}
\input{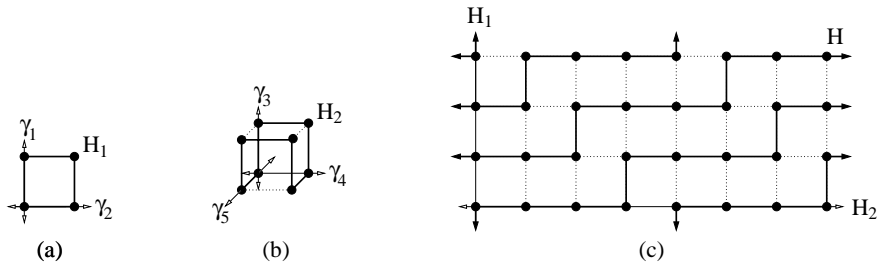}
\caption{Composing cycles. A Hamiltonian cycle $H_1$ in $\gamma_1
\times \gamma_2$ is combined with $H_2$ from $\gamma_3 \times
\gamma_4 \times \gamma_5$ to create a cycle $H$ in $\gamma_1
\times \gamma_2 \times \cdots \times \gamma_5$.  One edge is taken
from the direction of $H_1$, then three edges from $H_2$, and so
on; since both $H_1$ and $H_2$ have run length two, $H$ has also
run length two.} \label{comp}
\end{center}
\end{figure}

Our primary method for finding Hamiltonian cycles with high run
length consists of combining two or more Hamiltonian cycles to
construct another Hamiltonian cycle in their product; this method
is modeled after those in \cite{optrun} and \cite{bent}.  Say the
grid graph $G_1$ is the Cartesian product of cycles
\mbox{$\gamma_1 \times \gamma_2 \times \cdots \times \gamma_m$}
and $G_2$ is $\gamma_{m+1} \times \gamma_{m+2} \times \cdots
\times \gamma_k$, and let $H_1$ and $H_2$ be Hamiltonian cycles in
$G_1$ and $G_2$ respectively. Then any Hamiltonian cycle $H$ in
$H_1 \times H_2$ is also a Hamiltonian cycle in $G=G_1 \times
G_2$. We could consider $H$ a cycle in the $k$-dimensional torus
$\gamma_1 \times \gamma_2 \times \cdots \times \gamma_k$ or in the
two-dimensional torus $H_1 \times H_2$, as demonstrated in
Figure~\ref{comp}.  If each $\gamma_i$ in Figure~\ref{comp} is a
cycle on two vertices, and if $H_1$ is a cycle in $G_1 = \gamma_1
\times \gamma_2$ and $H_2$ is a cycle in $\gamma_3 \times \gamma_4
\times \gamma_5$, then an $H$ which is Hamiltonian in $H_1 \times
H_2$ is Hamiltonian in $G_1 \times G_2$.  We can even determine
the run length of $H$ from this representation. Note that $H$
takes three edges from the $H_2$ dimension followed by one from
$H_1$, and then repeats; if $rl(H_2) \geq 3$ and $rl(H_1) \geq 1$,
then any adjacent four edges come from different dimensions.
However, since $rl(H_2)=2$, any three adjacent edges from $H_2$
may repeat a dimension, and so the run length may be no greater
than two. We do know that any two adjacent edges are from
different dimensions and thus that $rl(H)$ is at least 2.  Note
that, to avoid difficulty, we can construct such an $H$ by
assuming $H_1$ and $H_2$ are directed.  In order to maximize
$rl(H)$, we repeat an edge from the same dimension of $H_i$ as
seldom as possible.

This construction suggests the following theorem, couched in terms
of $mrl(G)$, the maximum run length of any Hamiltonian cycle $H$
in $G$.  If $G=\gamma_1 \times \gamma_2 \times \cdots \times
\gamma_k$, then $|G|$, the number of vertices in the grid graph
$G$, is given by $\prod ^k _{i=1} |\gamma_i|$ .

\begin{thm}\label{st1}
Let $G_1$ and $G_2$ be grid graphs.  If there exist $s_1$ and
$s_2$ such that
\begin{enumerate}
\item$\gcd(|G_1|,s_1)=1,~\gcd(|G_2|,s_2)=1,$\label{cond1} and
\item$\gcd(|G_1|,|G_2|)=s_1+s_2$,\label{cond2}
\end{enumerate}
then $\mathrm{mrl}(G_1 \times G_2) \geq
\mathrm{mrl}(G_1)+\left\lfloor \frac{s_2}{s_1}
\mathrm{mrl}(G_1)\right\rfloor.$
\end{thm}

\begin{proof}
Choose Hamiltonian cycles $H_1$ and $H_2$ in $G_1$ and $G_2$,
respectively, with $rl(H_1)=mrl(G_1)$ and $rl(H_2)=mrl(G_2)$. Then
we construct a Hamiltonian cycle $H$ in $H_1 \times H_2$ with
greatest possible run length. In~\cite{ham1}, Trotter and
Erd\H{o}s show that a Hamiltonian cycle $H$ exists in the directed
Cartesian product $H_1 \times H_2$ exactly when there are integers
$s_1$ and $s_2$ satisfying Conditions~\ref{cond1} and~\ref{cond2}
above.

Curran and Witte explain in detail in \cite{ham2} that
Condition~\ref{cond2} ensures there is a cover of cycles where
each point is contained in exactly one cycle exactly once.  Any
$s_1 + s_2$ consecutive edges contain $s_1$ from the dimension
$H_1$ and $s_2$ from the dimension $H_2$.  In fact, any such
subsequence of $s_1 + s_2$ consecutive edges uniquely determines
the rest of the cycle or cover of cycles. Condition~\ref{cond1},
on the other hand, ensures that $H$ is a single Hamiltonian cycle
rather than a cover of cycles since $\gcd(|G_1|,s_1) \cdot
\gcd(|G_2|,s_2)$ gives the number of disjoint cycles in any such
cover.

Subject to these conditions, we describe a construction to
maximize $rl(H)$ based on the methods of Goddyn et
al.~\cite{optrun}. First, we pick $s_1$ and $s_2$ so that
$s_1/s_2$ is as close as possible to $rl(H_1)/rl(H_2)$.  We then
construct the cycle $H$ by adding edges in steps so that, at the
time of the $j^{th}$step, we have added $x_{1j}$ from the
direction of $H_1$ and $x_{2j}$ from the direction of $H_2$.  At
any point, the ratio $x_{1j}/x_{2j}$ is as close as possible to
$s_1/s_2$ without exceeding it.  That is, if we start with an edge
from the direction of $H_2$, we then add as many as possible from
the direction of $H_1$, i.e., we add $\lfloor s_1/s_2 \rfloor$
edges, so $(x_{11},x_{21})=(\lfloor s_1/s_2 \rfloor,1).$  Then we
add another from $H_2$ followed by the appropriate number from
$H_1$, so $(x_{12},x_{22})=(\lfloor 2s_1/s_2 \rfloor,2)$, and so
on. Ultimately, our first $s_1 + s_2$ edges will contain $s_1$
from $H_1$ and $s_2$ from $H_2$. Any subpath of these $s_1 + s_2$
edges containing exactly $x_{1j}=rl(H)$ edges from $H_1$ will
contain at least $x_{2i} \geq x_{1i}s_2/s_1 \geq \left\lfloor
\frac{s_2}{s_1}rl(H_1) \right\rfloor$ edges from $H_2$.  That is,
\[rl(H) \geq rl(H_1) + \left\lfloor \frac{s_2}{s_1} rl(H_1)
\right\rfloor.\] Note that interchanging $H_1$ and $H_2$ may give
a better bound.
\end{proof}

\begin{cor}\label{basiccor}
If there exist $s_1$ and $s_2$ satisfying the conditions of
Theorem~\ref{st1}, and if
\[\frac{s_1}{s_2}=\frac{\mathrm{mrl}(G_1)}{\mathrm{mrl}(G_2)},\]
\[then~\mathrm{mrl}(G_1 \times G_2) \geq \mathrm{mrl}(G_1) +
\mathrm{mrl}(G_2).\]
\end{cor}

As an application of Theorem~\ref{st1} and its
Corollary~\ref{basiccor}, we have the following results.

\begin{ex}\label{rlk}
If $G=\gamma_1 \times \gamma_2 \times \cdots \times \gamma_k$ is
directed, then $mrl(G) < k$.
\end{ex}
\begin{proof}
Recall that the method introduced in Theorem~\ref{st1} treats
$H_1$ and $H_2$ as if they are directed. The theorem of Trotter
and Erd\H{o}s gives necessary and sufficient conditions for the
existence of a Hamiltonian cycle in the product of exactly two
directed cycles; thus, the conditions imposed are sufficient for
the existence of cycles in $G$ if it is undirected but are also
necessary if it is considered as the product of directed cycles
$H_1 \times H_2$. In our case, any cycle in $G$ can be decomposed
into a cycle in $H_1 \times H_2$ for some directed $H_1$ and
$H_2$. That is, if any directed graph $\gamma_1 \times \gamma_2
\times \cdots \times \gamma_i$ has maximum run length equal to
$i$, then, for some ordering of the edges, there is a cycle in the
graph which consists of the repeated subsequence of one edge from
$\gamma_1$, then one from $\gamma_2$, $\dots$, one from
$\gamma_i$, another one from $\gamma_1$, etc.  This implies first
that $|G_1|=|\gamma_1||\gamma_2|\cdots|\gamma_{k-1}|$ is divisible
by $i=k-1$.  It also means, however, that some cycle in the $G$
given above can be decomposed into a cycle in $G_1 \times
\gamma_k$ where $s_1=mrl(G_1)=k-1$. Then we know that
$gcd(|G_1|,k-1)=1$, which is a contradiction.
\end{proof}

We may also make Corollary~\ref{basiccor} more specific.

\begin{cor}\label{adjustcor}
If $s_1$ and $s_2$ satisfy the conditions of Theorem~\ref{st1} and
if
\[|s_2 \cdot \mathrm{rl}(H_1) - s_1 \cdot \mathrm{rl}(H_2)| \leq
\mathrm{max}(s_1,s_2),\]then\[\mathrm{mrl}(G_1 \times G_2) \geq
\mathrm{mrl}(G_1) + \mathrm{mrl}(G_2) -1.\]
\end{cor}

\begin{proof}
Say $s_1 \geq s_2$. Then:\[\left|\frac{s_2}{s_1} \mathrm{rl}(H_1)
- \mathrm{rl(H_2)}\right| \leq 1\]so\[\left\lfloor
\frac{s_2}{s_1}\mathrm{rl}(H_1)\right\rfloor - \mathrm{rl}(H_2)
\leq 1.\] The corollary follows.
\end{proof}

Corollary~\ref{adjustcor} suggests another explicit application.

\begin{ex}\label{exp}
Let $p$ be a prime and let $G=\gamma_1 \times \gamma_2 \times
\cdots \times \gamma_{3n}$ with $|G|=p^r$ for some positive
integer $r$. Then $mrl(G) \geq n+1$.
\end{ex}
\begin{proof}
The base case was proved by Ruskey and Sawada~\cite{bent}, whose
results imply $mrl(\gamma_1\times\gamma_2\times\gamma_3)\geq 2$ if
$|\gamma_1||\gamma_2||\gamma_3|=p^r$.

For $n>1$, consider $j=\lceil\frac{n}{2}\rceil$ and let
$G_1=\gamma_1 \times \gamma_2 \times \cdots\times \gamma_{3j}$ and
$G_2=\gamma_{3j+1} \times \gamma_{3j+2} \times \cdots
\times\gamma_{3n}$ where $|G_1|=p^{r_1}$ and $|G_2|=p^{r_2}$ for
$p$ a prime and $r_1,~r_2$ integers at least 3.  Take Hamiltonian
cycles $H_1$ in $G_1$ and $H_2$ in $G_2$. Then $G=H_1 \times H_2$
and $\gcd(|H_1|,|H_2|)=p^r_0$ where $r_0=min(r_1,r_2)$.  Assume
inductively that $rl(H_1)=j+1$ and $rl(H_2)=n-j+1$. Then set
\[s_1=\left\lfloor\frac{\mathrm{rl}(H_1)}
{\mathrm{rl}(H_1)+\mathrm{rl}(H_2)}\mathrm{gcd}(|H_1|,|H_2|)
\right\rfloor+i=\left\lfloor\frac{j+1}{n+2}p^{r_0} \right\rfloor
+i\] with $i = 0$ or $1$ so $s_1$ is not divisible by $p$.  Also
set $s_2=p^{r_0}-s_1$ which is thus also not divisible by $p$.
Then the conditions of Theorem~\ref{st1} are satisfied, and
\[s_2\cdot \mathrm{rl}(H_1)-s_1\cdot \mathrm{rl}(H_2)|=|s_s(j+1)-s_1(n-j+1)| \leq 2
\leq s_1,\] so we may apply Corollary~\ref{adjustcor} to obtain
$mrl(H_1\times H_2)\geq n+1$.
\end{proof}

\begin{ex}\label{exp2}
Let $p$ be a prime and $G=\gamma_1 \times \gamma_2 \times \cdots
\times \gamma_{2n}$ where $|\gamma_i|=2^{q_i}p^{r_i}$ for integers
$q_i \geq 1$ and $r_i \geq 0$. Then $mrl(G) \geq n+1$.
\end{ex}
\begin{proof}
The method is very similar to that of Application~\ref{exp}.
Again, Ruskey and Sawada's work~\cite{bent} gives the base case
$mrl(\gamma_1\times\gamma_2)=2$ for $|\gamma_1|$ and $|\gamma_2|$
even.

Set $j=\lceil\frac{n}{2}\rceil$ and take Hamiltonian cycles $H_1$
in $G_1=\gamma_1 \times \gamma_2\times\dots\times\gamma_{2j}$ and
$H_2$ in $G_2=\gamma_{2j+1}\times\gamma_{2j+2}\times\dots
\times\gamma_{2n}$ so $G=H_1\times H_2$ and
$gcd(|H_1|,|H_2|)=2^{q_0}p^{r_0}$ for some integers $q_0\geq 1$
and $r_0\geq 0$.  Assuming inductively that $rl(H_1)=j+1$ and
$rl(H_2)=n-j+1$, we emulate Application~\ref{exp} and set
\[s_1=\left\lfloor\frac{j+1}{n+2}2^{q_0} p^{r_0}\right\rfloor+i\] for
$i=0,1,2,$ or 3 chosen so $s_1$ not divisible by 2 or $p$. Then
$s_2=2^{q_0}p^{r_0}-s_1$ is also not divisible by either. The
conditions of Theorem~\ref{st1} are met, and
\[|s_2(j+1)-s_1(n-j+1)|\leq 4\leq s_1,\] so
Corollary~\ref{adjustcor} applies and $rl(H_1 \times H_2) \geq
n+1$.
\end{proof}

\section{Adding extra dimensions} \label{secdim}

We might ask, though, if there are easier ways of bounding
$mrl(G)$. So far, we have primarily broken grid graphs into
constituents $G_1 \times G_2$ where each constituent has maximum
run length at least two, and where the run length of the whole is
greater than that of either constituent.  It is also possible,
however, to ask if the Cartesian product of $G$ and a
$k$-dimensional cycle has maximum run length at least that of $G$.
A method suggested by Ruskey and Sawada~\cite{bent} is useful for
one type of case.

\begin{thm}\label{extradim1}
If $G=\gamma_1 \times \gamma_2 \times \cdots \times \gamma_k$,
$|G|$ is even, and $\gamma_{k+1}$ is a cycle, then $mrl(G \times
\gamma_{k+1}) \geq mrl(G)$.
\end{thm}

\begin{proof}
Without loss of generality, say $|\gamma_1|=2n$.  Then consider
the embedding of $H$ in $G$ as a graph of connected lattice points
\textbf{G} in $\mathbb{R}$, which is also a graph in $\mathbb{R}^k
\times \{0\} \subset \mathbb{R}^{k+1}$.  We wish to stack
$|\gamma_{k+1}|$ copies of \textbf{G} on top of one another in
$\mathbb{R}^{k+1}$ with every other one reversed. Reflect
\textbf{G} exactly $i$ times in the line $x_1=n+1/2$ to create
$\textbf{G}_i$ and then set it in $\mathbb{R} \times \{ i \}$. Do
this until there are $|\gamma_{k+1}|$ copies of \textbf{G} fixed
immediately above each other with every other one reversed.

\begin{figure}[h]
\begin{center}
\input{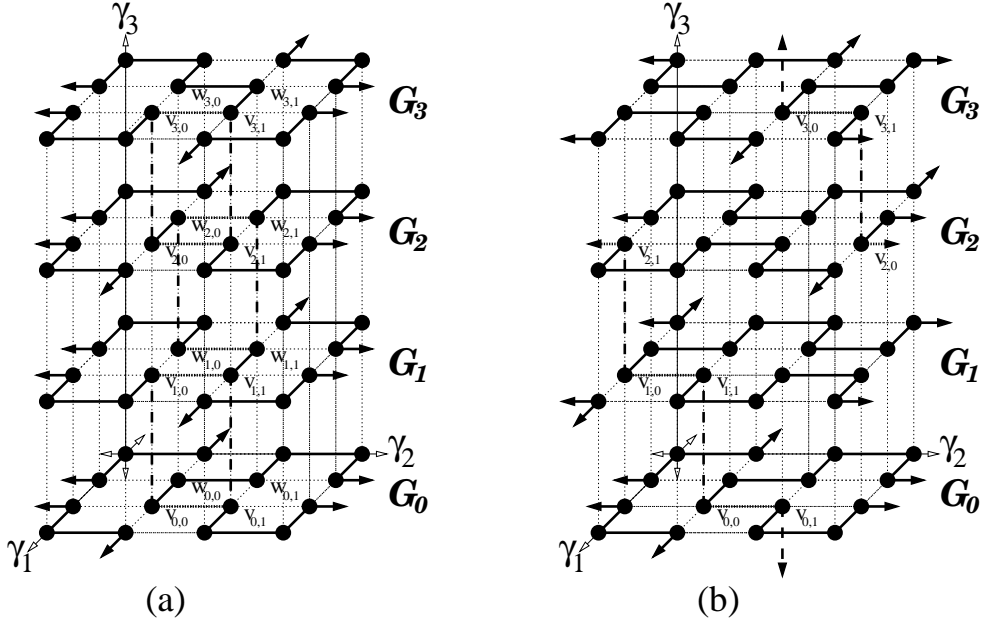}
\caption{Two methods for extending cycles without decreasing run
length are (a) extending the two-dimensional cycle $G$ from
$\gamma_1 \times \gamma_2$ into three dimensions where
$|\gamma_3|$ is even and (b) extending the two-dimensional cycle
$G$ from $\gamma_1 \times \gamma_2$ to a three-dimensional one
where $|\gamma_3|$ is a multiple of $|\gamma_2|$. In both cases,
the solid lines represent the original cycle $H$, the dashed lines
mark omitted edges, and the dotted lines mark added edges.}
\label{extradim}
\end{center}
\end{figure}

Pick any two edges $[v_{i,0}, v_{i,1}]$ and $[w_{i,0},w_{i,1}]$
from the $\gamma_1$ direction located in the middle of
$\textbf{G}_0$. Then the point at $v_{i,0}$ is $\langle n,x_2,
x_3, \cdots, x_k, i\rangle$ and \mbox{$v_{i,1}=v_{i,0}+\langle
1,0,0,\cdots,0\rangle$,} and, similarly, $w_{i,0}=\langle n,y_2,
y_3, \cdots, y_k, i\rangle$ and $w_{i,1}=w_{i,0}+\langle
1,0,0,\cdots,0\rangle.$ We assume that $x_i \neq y_i$. Note that
these edges are preserved by the reflections, so $[v_{i,0},
v_{i,1}]$ is directly below $[v_{i+1,0}, v_{i+1,1}]$, and
$[w_{i,0},w_{i,1}]$ is also below $[w_{i+1,0},w_{i+1,1}]$. Then
remove these four edges and insert $[v_{i,0},
v_{i+1,0}],~[v_{i,1}, v_{i+1,1}],~[w_{i,0},w_{i+1,0}],$ and
$[w_{i,1},w_{i+1,1}]$.  This creates a cycle in $\bigcup
\textbf{G}_i$, which in turn induces a cycle $H'$ in $G\times
\gamma_{k+1}$. Where $rl(H)=mrl(G)$, any subsequence of $rl(H)$
adjacent edges in $H$ with no more than one edge from any one
dimension now corresponds to a sequence of $rl(H)$ adjacent edges
in $H'$ with still no more than one edge from any one direction.
Thus, $mrl(G\times\gamma_{k+1})\geq mrl(G)$ See
Figure~\ref{extradim}a.
\end{proof}

As a continuation from Application~\ref{exp2}, we have our next
result.

\begin{ex}
Let $G=\gamma_1 \times \gamma_2 \times \cdots \times \gamma_k$ and
$|\gamma_i|=2^{q_i}p^{r_i}$, $q_i \geq 1, r_i \geq 0.$ Then
$mrl(G) \geq \lfloor k/2 \rfloor+1$ for all $k$.
\end{ex}
\begin{proof}
If $k$ is even, we may apply Corollary~\ref{adjustcor}; if $k$ is
odd, $k-1$ is even, and we may apply Theorem~\ref{extradim1} and
then Corollary~\ref{adjustcor} in succession.
\end{proof}

A variation of this method is still available even if all the
$|\gamma_i|$ are odd.

\begin{thm}\label{extradim2}
If $G=\gamma_1 \times \gamma_2 \times\cdots \times \gamma_k$ and
$|\gamma_{k+1}|$ is a multiple of $|\gamma_j|$ for some $j \leq
k$, then $mrl(G \times \gamma_{k+1}) \geq mrl(G)$.
\end{thm}

\begin{proof}
The proof is very similar to that of Theorem~\ref{extradim2}.  The
graph $G$ embeds nicely as $\mathcal{G}_i$ in $\mathbb{R}^k \times
\{ i \}$; however, instead of reflecting $\mathcal{G}$, we
translate it by $i$ times the $j^{th}$ unit vector.  Then we
remove some $[v_{1,0},v_{1,1}]$ in $\mathcal{G}_1$ and all its
images $[v_{i,0},v_{i,1}]$ in $\mathcal{G}_i$.  We attach
$[v_{i,1},v_{i+1,0}]$ to get a cycle in $\bigcup \mathcal{G}_i$
with the same run length as $G$.  (Note that we must make sure
that $v_{|\gamma_j|+1} = v_0$.)  See Figure~\ref{extradim}b.
\end{proof}

Again, Theorem~\ref{extradim2} allows us to generalize
Application~\ref{exp} as follows.

\begin{ex}
If $G=\gamma_1\times\gamma_2\times\cdots\times\gamma_k$ and
$|G|=p^r$ for positive integer $r$, then $mrl(G) \geq \lfloor k/3
\rfloor +1$ for all $k$.
\end{ex}
\begin{proof}
We already know that $mrl(G) \geq k/3 +1$ for all $k$ divisible by
three.  The rest follows by Theorem~\ref{extradim2}.
\end{proof}

\section{Conclusion and further work}

This paper sets lower bounds on the maximum run length of some
finite dimensional tori by giving methods for constructing
Hamiltonian cycles of known run length. These methods produce
Hamiltonian cycles in a torus $G=G_1 \times G_2$ from Hamiltonian
cycles $H$ in $G_1$ and $G_2$.

The methods do not give very close upper and lower limits on
$rl(H)$ and $mrl(G)$ and may be improved by future research.
Additionally, there are many interesting cases not addressed in
this paper, including: those cases where a path or a cover of
disjoint cycles, but not a Hamiltonian cycle, can be found; grid
graphs which are the cartesian products of paths or of paths and
cycles; and grid graphs where some dimensions are directed and
others are undirected. It would also be interesting to work out
exactly which portion of the many possible grid graphs are covered
by this work (since many of the results depend on the greatest
common divisors of the sizes of the factor graphs).  Additionally,
the methods of Theorem~\ref{st1} may be adjustable for any
$k$-dimensional grid graph, perhaps using an extension of the idea
of a diagonal (presented by Curran and Witte in \cite{ham2}).
Goddyn et al.~\cite{optrun} do present a version for a
three-dimensional grid graph on eight vertices.  The results in
this paper also suggest a generalization of run length to
\emph{j-variance}, which we define to be the maximum $r$ such that
any sequence of $r$ consecutive edges in a cycle $H$ includes at
most $j$ edges from any one dimension (\emph{N.B.:} 1-variance is
the same as run length). Finally, a general but difficult
challenge is finding methods for determining or estimating the
maximum run length (or perhaps $j$-variance) of an arbitrary grid
graph.

\section{Acknowledgements}
The work for this paper was mostly completed at the University of
Minnesota, Duluth, funded by NSA grant MDA~904-02-1-0060 and NSF
grant DMS-0137611.  Many thanks are owed to the unflagging
interest and input of Joe Gallian, Geir Helleloid, and Phil
Matchett.

\end{document}